\documentclass{amsart}
\usepackage{amsmath}
\usepackage{amsfonts}
\usepackage{amssymb}
\usepackage{amsthm}
\usepackage{url}
\usepackage[all]{xy}
\usepackage{dsfont}
\usepackage{graphicx}
\usepackage{caption}
\usepackage{subcaption}
\usepackage{comment} 
\usepackage{stmaryrd}
\usepackage{hyperref}
\usepackage{todonotes}
\usepackage{color}
\usepackage{enumerate,tikz-cd}
\usepackage{fixltx2e}
\usepackage{MnSymbol}
\usepackage{comment}
\usepackage{mathtools}
\allowdisplaybreaks

\numberwithin{equation}{section}

\newtheorem{proposition}{Proposition}[section]
\newtheorem{lemma}[proposition]{Lemma}
\newtheorem{theorem}[proposition]{Theorem}

\theoremstyle{definition}

\newtheorem{definition}[proposition]{Definition}
\newtheorem{example}[proposition]{Example}

\DeclareMathOperator{\rk}{rk}

\renewcommand{\epsilon}{\varepsilon}

\renewcommand{\phi}{\varphi}

\newcommand\mi{^{-1}}

\newcommand\bsni{\bigskip\noindent}

\newcommand\bbc{\mathbb{C}}

\newcommand\bbp{\mathbb{P}}

\newcommand\cA{\mathcal{A}}

\newcommand\cO{\mathcal{O}}

\newcommand\cT{\mathcal{T}}

\newcommand\bigbbp{{\underline{\bbp}}}

\title[Deligne pairings for projective duals]{Deligne pairing isomorphisms for projective duals and higher-rank associated hypersurfaces}
\author{Rémi Reboulet}

\address{Rémi Reboulet, CNRS, Institut Camille Jordan, 21 Av. Claude Bernard, 69100 Villeurbanne, France}\email{reboulet@math.univ-lyon1.fr}

\begin{document}

\maketitle

\begin{abstract}
	Given a family $\pi:X\to S$ of complex projective varieties, polarised by a (sufficiently) relatively very ample $L$, we give an isomorphism between a Deligne pairing involving higher relative cotangent Chern classes, and the relative projective dual line $X^\vee(L)$. As a consequence, we obtain explicit Deligne pairing isomorphisms for all of the higher associated hypersurfaces of Gelfand--Kapranov--Zelevinsky, including mixed polarisations, as well as higher-rank versions of such hypersurfaces involving sufficiently positive vector bundles.
\end{abstract}

\tableofcontents

\section*{Introduction.}

This note aims at showing a simple way to generalise a class of results yielding isomorphisms between \textit{resultant lines} and \textit{Deligne pairings}. The first notable result of this type is due to Zhang \cite{zhang:height} (see also similar previous results in Arakelov geometry: \cite{arakelov:i,arakelov:ii,arakelov:iii,arakelov:iv,arakelov:v}), who constructed an isomorphism between the Chow line and the top-Deligne pairing $\langle L^{n+1}\rangle_{X/S}$ of a very ample line bundle on a family $\pi:X\to S$. Yhis result saw a slight generalisation in \cite{qtian}, with the full mixed case $\langle L_1,\dots,L_n\rangle_{X/S}$ treated in \cite{reb:4} (see also \cite{chen:moriwakiintersection} in the Arakelov setting).

Work of Kapadia \cite{kapadia,kapadia:thesis} then extended Zhang's result to a modification of the Chow line, the \textit{hyperdiscriminant} (or \textit{Hurwitz}) \textit{line}, which may be constructed using the first associated hypersurface described in \cite{book:gkz}. This, in view of Paul's \cite{paul:chow,paul:hyperdiscrimants} and Phong--Sturm's \cite{phongsturm:stability} results, was motivated by the study of K-stability. Analogous ideas were used in \cite{dr:3} and \cite{reb:4}.

We show in this note that one can describe arbitrary higher associated hypersurfaces in the sense of \cite{book:gkz}, as well as vector bundle resultants \cite{abolazsmith}, and many other constructions, in terms of Deligne pairings. All of those results may be unified simply by realising the (relative) projective dual of the embedding of $X$ by a single relatively very ample line bundle, as a Deligne pairing. Namely, we prove the following.

\bsni\textbf{Theorem A.} Let $\pi:X\to S$ be a family satisfying hypotheses $(H_n)$ as below, and let $L$ be a (sufficiently) relatively very ample line bundle on $X$. Let $\cA$ be the line defining the relative projective dual $X^\vee(L)$. Then, there is a canonical isomorphism
	\begin{equation}
		\cA\simeq\left\langle c_1(L)\cdot c_n(T_{X/S}^\vee \otimes L)\right\rangle_{X/S}.
	\end{equation}

\bigskip The proof follows the ideas of \cite{zhang:height,reb:4} and \cite{kapadia}, adequately modified to fit the construction: we evaluate an auxiliary pairing along two different sections, which define birational hypersurfaces, and use standard properties of the Deligne pairing. Applying Theorem A to the hyperplane bundles of projectivised vector bundles of the form $\bigoplus_j E_j$, with $E_j$ positive vector bundles on $X$, then yields a wealth of applications, described at the end of the article. We obtain, in particular, an expression for the higher associated hypersurfaces of \cite{book:gkz}
\begin{equation}
\cA_{n+1-r}=\langle c_1(L)^r \cdot c_{n+1-r}(T_{X/S}^\vee\otimes L)\rangle_{X/S}.
\end{equation}
as well as a higher-rank version of the main result of \cite{reb:4}: given very positive vector bundles $E_1,\dots,E_r$ on $X$, with degrees summing to $n+1$:
\begin{equation}
\cA = \langle c_{\rk(E_1)}(E_1)\cdot\dots\cdot c_{\rk(E_r)}(E_r)\rangle_{X/S}.
\end{equation}
We refer the reader to the last section for more precise statements.


\bsni\textbf{Acknowledgements.} The author thanks Ruadhaí Dervan for discussions related to the main results of this article.

\bsni\textbf{AI usage statement.} The entirety of this article was written by the author. It is based on an old idea of using Kapadia's arguments to prove the general isomorphism, and was considered in \cite{reb:4} (see Section 3.6 therein), but ultimately left for a later article. LLMs were used in verifying and clarifying calculations in the present article, as well as to clarify results in \cite{kapadia}, and auditing the article for correctness.

\section{Preliminaries.}

\subsection{Deligne pairings for Chern polynomials.}

Deligne pairings were constructed by Elkik \cite{elkik:1} (building on \cite{deligne:detcoh}); we give their construction using generators and relations as in \cite[II.2.4]{elkik:1} and \cite{zhang:height}. For the entirety of this section, we refer the reader to the aforementioned \cite{elkik:1,zhang:height}, as well as \cite{eriksson:freixas}. We recall afterwards how to construct Deligne pairings for more general Chern polynomials, as will be needed for us.

Our standing assumption on families $X\to S$ throughout the article will be the following:
\begin{definition}\label{def:hn}We say a morphism $\pi:X\to S$ satisfies $(H_n)$ if it is smooth, flat, projective, finitely presented, and of relative (pure) dimension $n$.
\end{definition}

We now fix $\pi:X\to S$ satisfying $(H_n)$ as above. We begin with the datum of line bundles $L_1,\dots,L_{n+1}$ on $X$. We then define the \textit{Deligne pairing}
$$\langle c_1(L_1)\cdot \ldots \cdot c_1(L_{n+1})\rangle_{X/S}$$
to be the line bundle on $S$ which has, as sections, formal expressions of the form
$$\langle s_1,\dots,s_{n+1}\rangle_{X/S},$$
with each $s_i\in H^0(X,L_i)$ and each collection satisfying that $\bigcap_{i=1}^{n+1} V(s_i)=\emptyset$, subject to all relations
$$\langle s_1,\dots,fs_{n+1},\dots,s_n\rangle_{X/S}=\sum_k a_k N_{D_k/S}(f)\langle s_1,\dots,s_{n+1}\rangle$$
where $i\in[n+1]$ is such that the divisor $\bigcap_{j\neq i}(s_j)=\sum a_k D_k$ is flat over $S$, and $f$ is a rational function on $X$. This turns out to quotient to a line bundle on the base $S$. We warn the reader that the notation for Chern classes may be confusing in case of repeated indices: $c_1(L^{\otimes k})=c_1(kL)$ denotes the genuine Chern class of $kL$, whereas $c_1(L)^{k}$ will denote the $k$-th product of $c_1(L)$, a degree $k$ Chern polynomial expression.

Let now $(E_i)_{i\in [r]}\to X$ be a collection of vector bundles with $\rk(E_i):=r_i+1$, and let $\bigbbp:=\bbp E_1\times_X \dots\times_X \bigbbp E_r$, with projections $p_i:\bigbbp\to \bbp E_i$, and set $H_i:= p_i^*\cO_{\bbp E_i}(1)$. Let $(k_i)_{i\in [r]}$ be nonnegative integers summing to $n+1$. One then defines Deligne pairings for Segre classes by
\begin{equation}
\langle s_{k_1}(E_1)\cdot\ldots\cdot s_{k_r}(E_r)\rangle_{X/S}:=\langle c_1(H_1)^{r_1+k_1}\cdot \ldots\cdot  c_1(H_r)^{r_r+k_r}\rangle_{\bigbbp/S}.
\end{equation}
We may have $k_j=0$ for some $j$, in which case $s_0(E_j)$ will simply denote omission in the Deligne pairing. It is then clear that we recover the usual construction for line bundles. We note, as observed in a more general setting in \cite[Remarque V.1.2(e)]{elkik:1}, that in particular if $r=n+1$,
\begin{equation}
	\langle s_{1}(E_1)\cdot\ldots\cdot s_{1}(E_{n+1})\rangle_{X/S}=\langle c_1(\det E_1)\cdot\ldots\cdot c_1(\det E_{d+1})\rangle_{X/S}.
\end{equation}

We now turn to Chern classes of the $E_i$, following Eriksson--Freixas \cite[Section 7.1.1]{eriksson:freixas} which coincides with Elkik's (by \cite[Corollary 7.11]{eriksson:freixas}). First set $c_0=s_0=1$, and recall the identities
$$c_k=\sum_{i=1}^k (-1)^{i+1}s_ic_{k-i}.$$
This allows us to define pairings for Chern classes inductively. Assuming the $k_j$'s to sum adequately, we set
\begin{equation}
\langle c_{k_1}(E_1)\cdot \ldots \cdot c_0(E_j)\cdot\ldots\cdot c_{k_r}(E_r)\rangle_{X/S}=\langle c_{k_1}(E_1)\cdot\ldots\cdot  s_0(E_j)\cdot \ldots\cdot c_{k_r}(E_r)\rangle_{X/S},
\end{equation}
\begin{equation}
\langle c_{k_1}(E_1)\cdot\ldots\cdot  c_k(E_j)\cdot \ldots\cdot c_{k_r}(E_r)\rangle_{X/S}=\sum_{p=1}^k(-1)^{p+1}\langle c_{k_1}(E_1)\cdot\ldots\cdot s_p(E_j)\cdot c_{k-p}(E_j)\cdot \ldots\cdot c_{k_r}(E_r)\rangle_{X/S}.
\end{equation}
It is then clear that, if all the $k_i=1$,
\begin{equation}
\langle c_{1}(E_1)\cdot\ldots\cdot c_1(E_{n+1})\rangle_{X/S}=\langle s_{1}(E_1)\cdot\ldots\cdot s_1(E_{n+1})\rangle_{X/S},
\end{equation}
which is further computed by determinants as above. The Deligne pairing is then extended, by linearity, to general Chern--Segre polynomials.

Finally, we note that we will use additive notation for the Deligne pairings themselves, as line bundles.

\subsection{Properties of Deligne pairings.}

We collect here a list of properties which will be used in the proof of the main theorem. Throughout, unless otherwise stated, $L$ will be used to denote line bundles, $E$ vector bundles, $P,Q$ homogeneous Chern polynomials, sometimes understood to have suitable degrees. The equal sign will always denote "there exists an isomorphism, compatible with the previous operations, which is natural in the sense of \cite{eriksson:freixas}". All of the following properties may be found in \cite{eriksson:freixas} and \cite{elkik:1}, to which we refer the reader.

We first have the following normalisations:
\begin{align}
	\langle 0\rangle_{X/S}&=\cO_S,\\
	\label{eq:normalisation} \langle c_1(L)\rangle_{S/S}&=L.
\end{align}
For integers $r,s$, $\deg P = \deg Q = n+1$,
\begin{equation}\label{eq:linearity}
	\langle rP+sQ\rangle_{X/S}=r\langle P\rangle_{X/S}+s\langle Q\rangle_{X/S}.
\end{equation}
For $\rk(E)=e$, $0\leq k\leq e$, $\deg P=n+1-k$,
\begin{equation}\label{eq:duality}
	\langle P\cdot c_k(E^\vee)\rangle_{X/S}=(-1)^k\langle P\cdot c_k(E)\rangle_{X/S}. 
\end{equation}
We note the following twisting formula:
\begin{equation}\label{eq:twisting}
	\langle P \cdot c_k(E\otimes L)\rangle_{X/S}=\sum_{i=0}^k \binom{e-i}{k-i}\langle P,c_i(E),c_1(L)^{k-i}\rangle_{X/S}.
\end{equation}
Suppose $\pi:X\to S$ satisfies $(H_n)$, and $\mu:X'\to X$ satisfies $(H_m)$ with $\pi\circ\mu:=\pi'$ satisfying $(H_{n+m})$. Let $P$ be a homogeneous Chern polynomial on $X$, $Q$ on $X'$, such that $\deg P + \deg Q = n+m+1$. Then,
\begin{align}
	\label{eq:proj1}\langle \mu^*P\cdot Q\rangle_{X'/S}&=\langle P\cdot c_1(\langle Q\rangle_{X'/X})\rangle_{X/S},&\quad \deg Q = m+1;\\
	\label{eq:proj2}\langle \mu^*P\cdot Q\rangle_{X'/S}&=\nu \cdot \langle P\rangle_{X/S},&\quad \deg Q=m;\\
	\label{eq:proj3}\langle \mu^*P\cdot Q\rangle_{X'/S}&=\cO_S,&\quad \deg Q<m,
\end{align}
where 
\begin{equation} 
\nu=\int_{X'_x}Q|_{X'_x}.
\end{equation}
Suppose $0\to E' \to E \to E''\to 0$ is exact and $P$ is homogeneous of degree $n+1-k$, then
\begin{equation}\label{eq:whitney}
\langle P\cdot c_k(E)\rangle_{X/S}=\sum_{i=0}^k \langle P\cdot c_i(E')\cdot c_{k-i}(E'')\rangle_{X/S}.
\end{equation}
In particular, for $\deg P$ adequately chosen,
\begin{align}
\langle P\cdot c_1(E)\rangle_{X/S}&=\langle P\cdot c_1(\det E)\rangle_{X/S},\\
\langle P\cdot c_0(E)\rangle_{X/S}&=\langle P\rangle_{X/S},
\end{align}
and if $k> \rk E$,
\begin{equation}
\langle P\cdot c_k(E)\rangle=\cO_S.
\end{equation}
Suppose $s$ is a regular section of $E$, $\rk E=e+1$, with vanishing locus $j:Z:=V(s)\hookrightarrow X$ flat over $S$ of relative dimension $n-e-1$. Then
\begin{equation}\label{eq:sect_restr}
\langle P \cdot c_{e+1}(E)\rangle_{X/S}=\langle j^*P\rangle_{Z/S}.
\end{equation}
In particular, for a flat relative effective Cartier $D$, applying \eqref{eq:sect_restr} to the canonical section $\cO_X(D)$ yields:
\begin{equation}\label{eq:relcartier}
	\langle P\cdot c_1(\cO_X(D))\rangle_{X/S}=\langle P|_D\rangle_{D/S}.
\end{equation}
Suppose $X,X'\to S$ both satisfy $(H_n)$ over $S$. Let $\mu:X'\to X$ be an isomorphism over a quasi-compact open of $X$, dense in every fibre (we will call such maps \textit{relative birational maps}, and say $X$ and $X'$ are \textit{relatively birational} over $S$). If $P$ is a Chern polynomial of degree $n+1$ on $X$,
\begin{equation}\label{eq:relbir}
\langle P\rangle_{X/S}=\langle \mu^*P\rangle_{X'/S}.
\end{equation}
Finally, we note a useful comparison for Segre classes, which is a particular case of the definitions in the previous section, is the following: let $E\to X$ be a vector bundle of rank $e+1$, and $k\geq 0$. Write $p:\bbp E \to X$. Let $P$ be a Chern polynomial of adequate degree on $X$. Then,
\begin{equation}\label{eq:delignesegre}
\langle p^*P \cdot c_1(\cO_{\bbp E}(1))^{e+k}\rangle_{\bbp E/S} = \langle P\cdot s_k(E)\rangle_{X/S}.
\end{equation}

\section{Deligne pairing isomorphism for the projective dual.}

\subsection{Relative dual line and incidence constructions.}

Let $\pi:X\to S$ satisfy $(H_n)$. Let $L$ be a relatively very ample line bundle on $X$, such that $\pi_*L$ is locally free of rank $N+1$. Let $p:\bbp:=\bbp(\pi_*L)\to S$. We denote by $X^\vee(L)$ the relative projective dual of $X$ by the relative embedding given by $L$. In other words, given the incidence
$$Z:=\{(x,[s]),\,s(x)=0,\,d_x s=0\}\subset X\times \bbp,$$
then $X^\vee$ is the image of $Z$ by the projection to $\bbp$, with its reduced scheme structure. We will denote by $X_s^\vee(L_s)$ its fibre, which is just the projective dual of $X_s$. We first note the following:
\begin{lemma}
For all $k$ large enough, uniformly in $s$, $X_s^\vee(kL_s)$ is a hypersurface.
\end{lemma}
\begin{proof}
By \cite[Chapter 1, Corollary 1.2]{book:gkz}, it suffices to show that the image of each $X_s$ by $M:=kL$, $k\geq 0$, does not contain a projective line. Suppose it is the case, and let $C_s$ be the preimage of such a line in $X_s$, then $M_s\cdot C_s = 1$. But since $L_s$ is ample, we have $L_s\cdot C_s\geq 1$, hence $M_s\cdot C_s = k(L_s\cdot C_s)\geq k$. Thus the result holds for any $k\geq 2$.
\end{proof}

In the case that $X^\vee(L)$ is a hypersurface in each fibre, it is automatically a relative effective Cartier divisor on $\bbp/S$. There then exists $d$ such that 
$$\cO_{\bbp_s}(X_s^\vee(L_s))\simeq \cO_\bbp(d).$$
We then define
\begin{equation}\label{eq:defmainline}
\cA:=p_*\left(\cO_\bbp(X^\vee(L)-d)\right)
\end{equation}
We will call $\cA$ the \textit{relative dual line}. By the projection formula, one then has
\begin{equation}\label{eq:defmainline2}
	\cO_\bbp(X^\vee(L))\simeq p^*\cA + \cO_\bbp(d).
\end{equation}
We now give a different incidence construction closely related to $X^\vee(L)$. Consider the jet bundle $J^1_{X/S}(L)$, and note that one has a canonical section given by jet evaluation,
$$\sigma\in H^0(X\times \bbp,\pi_X^* J^1_{X/S}(L)\otimes \pi_\bbp^* \cO_\bbp(1))$$
with $\left(\sigma(x,[t]):\lambda t\mapsto j_x^1(\lambda t)\right)\in (J^1_{X/S})_x\otimes (\bbc t)^\vee$. By very ampleness, the evaluation morphism at each $x\in X$ is surjective, and the vanishing locus $V(\sigma)\subset X\times\bbp$ then has as fibre over $x$ its projectivised kernel $\bbp \ker\left(ev^1_x:H^0(X_s,L_s)\to (J_{X/S}^1)_s\right)$. We also will use the following:
\begin{proposition}\label{prop:relbir}
	The projection $\pi_\bbp|_{V(\sigma)}: V(\sigma)\to X^\vee(L)$ is relatively birational over $S$.
\end{proposition}
\begin{proof}
Let $U$ be the smooth locus of $X^\vee(L)$, which is dense in every fibre of $X^\vee(L)$ by \cite[Chapter 1, 1]{book:gkz}. Then by \cite[Chapter 1, Proposition 3.2(b)]{book:gkz}, $\pi_\bbp|_{V(\sigma)}\mi(U)$ is a fibrewise isomorphism, hence an isomorphism, concluding the proof.
\end{proof}
By \cite[Chapter 3, Example 3.12]{book:gkz}, one then has that
\begin{equation}\label{eq:degree}
d = \int_{X_s} c_n(J_{X_s}^1(L_s)) = \sum_{p=0}^n (n+1-p)\int_{X_s}c_p(T^\vee_{X_s})\cdot c_1(L_s)^{n-p}.
\end{equation}


\subsection{Proof of the main result.}

We now prove the following.
\begin{theorem}\label{thm:main}
	Suppose $(\pi:X\to S,L)$ are as in the previous subsection, and that $X^\vee(L)$ is a relative effective Cartier divisor. We then have that
	\begin{equation}\label{eq:main}
		\cA=\left\langle c_1(L)\cdot c_n(T_{X/S}^\vee \otimes L)\right\rangle_{X/S}.
	\end{equation}
\end{theorem}
\begin{proof}
Consider the auxiliary pairing
\begin{equation}
	\cT:=\left\langle c_1(\pi_\bbp^*\cO_\bbp(1))^N\right\rangle_{V(\sigma)/S}.
\end{equation}
We will first evaluate it along the relative birational morphism defined above, then along $\sigma$. Using thus \eqref{eq:relbir} (and Proposition \ref{prop:relbir}), then \eqref{eq:relcartier}, we obtain
\begin{equation}
	\cT=\langle c_1(\cO_\bbp(1))^N\rangle_{X^\vee(L)/S}=\langle c_1(\cO_\bbp(1))^N\cdot c_1(\cO_\bbp(X^\vee(L)))\rangle_{\bbp/S}. 
\end{equation}
By \eqref{eq:defmainline2}, linearity, and finally \eqref{eq:proj2} together with $\int_{\bbp_s} c_1(\cO_\bbp(1))^N=1$, we find
\begin{align}
\cT&=\langle c_1(\cO_\bbp(1))^N\cdot c_1(p^*\cA)\rangle_{\bbp/S}+d\langle c_1(\cO_\bbp(1))^{N+1}\rangle_{\bbp/S}\\
&=\langle c_1(\cA)\rangle_{S/S} + d\langle c_1(\cO_\bbp(1))^{N+1}\rangle_{\bbp/S}.
\end{align}
By \eqref{eq:normalisation} this is finally
\begin{equation}\label{eq:comp1}
\cT=\cA + d\langle c_1(\cO_\bbp(1))^{N+1}\rangle_{\bbp/S}.
\end{equation}
Using now \eqref{eq:sect_restr} applied to $\sigma$, and \eqref{eq:twisting},
\begin{align}
\cT&=\langle c_1(\pi_\bbp^*\cO_\bbp(1))^N\cdot c_{n+1}(\pi_X^*J^1_{X/S}(L) \otimes\pi_\bbp^*\cO_\bbp(1))\rangle_{X\times\bbp/S}\\
&=\sum_{i=0}^{n+1}\langle c_{n+1-k}(\pi_X^*J^1_{X/S}(L))\cdot c_1(\pi_\bbp^*\cO_\bbp(1))^{N+k}\rangle_{X\times\bbp/S}.
\end{align}
By \eqref{eq:proj3}, only the terms with $k<2$ are nonvanishing; the projection formula \eqref{eq:proj2} applied to the two remaining terms then yields
\begin{equation}\label{eq:comp2}
\cT=\langle c_{n+1}(J^1_{X/S}(L))\rangle_{X/S} + \int_{X_s}c_n(J^1_{X_s}(L))\langle c_1(\cO_\bbp(1))^{N+1}\rangle_{\bbp/S}.
\end{equation}
Using \eqref{eq:degree}, comparing \eqref{eq:comp1} and \eqref{eq:comp2} yields
\begin{equation}
\cA = \langle c_{n+1}(J^1_{X/S}(L))\rangle_{X/S}.
\end{equation}
To conclude, we finally apply \eqref{eq:whitney} to the exact sequence
$$0\to T_{X/S}^\vee \otimes L \to J^1_{X/S}(L)\to L\to 0,$$
where only the term with $c_1(L)$ is nonvanishing. This concludes the proof.
\end{proof}

\section{Important special cases.}


We prove the following as a consequence of Theorem \ref{thm:main}. As we will see, its corollaries will recover more familiar objects.

\begin{theorem}\label{thm:vb}
Let $\pi:X\to S$ satisfy $(H_n)$. Let $E$ be a vector bundle of rank $e$ on $X$. Write $q:Y:=\bbp E\to X$, $M:=\cO_{\bbp E}(1)$, and assume $Y^\vee(M)$ is relative effective. Let $\cA$ be the relative dual line $Y^\vee(M)$. Then,
$$\cA=\left\langle c_{e}(E)\cdot \sum_{p=0}^{n+1-e}c_p(T^\vee_{X/S})\cdot s_{n+1-e-p}(E)\right\rangle_{X/S}.$$
\end{theorem}
\begin{proof}
By Theorem \ref{thm:main}, we have
\begin{equation}
\cA = \langle c_1(M)\cdot c_{n+e-1}(T^\vee_{Y/S}\otimes M)\rangle_{Y/S}.
\end{equation}
Applying \eqref{eq:whitney} to the (twisted) cotangent sequence
$$0\to q^*T_{X/S}^\vee\otimes M\to T_{Y/S}^\vee\otimes M\to T_{Y/X}^\vee\otimes M\to 0,$$
and using that only the term with $i=n$ is nonvanishing, yields
\begin{equation}
\cA = \langle c_1(M)\cdot c_n(T^\vee_{X/S}\otimes M)\cdot c_{e-1}(T^\vee_{X/Y}\otimes M)\rangle_{Y/S},
\end{equation}
while \eqref{eq:whitney} now applied to the relative Euler sequence
$$0\to T^\vee_{Y/X}\otimes M\to q^*E\to M\to 0$$
with only the case $i=e-1$ nonvanishing, yields
\begin{equation}
\cA = \langle c_e(q^*E)\cdot c_n(q^*T^\vee_{X/S}\otimes M)\rangle_{Y/S}.
\end{equation}
Expanding via \eqref{eq:twisting},
\begin{align}
\cA&=\left\langle  c_e(q^*E)\cdot \left(\sum_{p=0}^n c_p(q^*T^\vee_{X/S})\cdot c_1(M)^{n-p}\right)\right\rangle_{Y/S},
\end{align}
and the result follows from \eqref{eq:delignesegre}.
\end{proof}

\begin{example}Let $E=\bigoplus_{j=1}^r E_j$ with the $E_1,\dots,E_r$ sufficiently very ample on $X$ and with ranks summing to $n+1$. Using the Whitney formula \eqref{eq:whitney} applied to the direct sum, one obtains from Theorem \ref{thm:vb} the higher-rank case of the results of \cite{reb:4}, alluded to in the last section thereof:
\begin{equation}
\cA = \langle c_{\rk(E_1)}(E_1)\cdot\dots\cdot c_{\rk(E_r)}(E_r)\rangle_{X/S}.
\end{equation}
In particular, for sufficiently very ample line bundles, one recovers the isomorphism part of \cite[Theorem B]{reb:4}, itself generalising \cite{zhang:height}. In the case where $n=1$, the vector bundle formula recovers the $IC_2(E)=\langle c_2(E)\rangle_{X/S}$ pairing of \cite{deligne:detcoh}, \cite[Definition 3.4]{eriksson:ic2}. The vector bundle resultant, i.e.\ the left-hand side above, also appears in \cite[Chapter 3, Example 3.9]{book:gkz}.
\end{example}

\begin{example}For $L$ sufficiently very ample and $E=L^{\oplus r}$ with $r\leq n$, $Y^\vee(M)$ is the (relative) $n+1-r$-th associated hypersurface of $(X,L)$ in the sense of \cite[Chapter 3, 2.E]{book:gkz}, that is, the set of codimension $r$ linear subspaces $W$ of $\bbp^N$ such that, for some $x\in X\cap W$,
$$\dim ET_x X\cap W \geq n+1-r.$$
Applying Theorem \ref{thm:vb} using $s_k(E) = \binom{r+k-1}{k}c_1(L)^k$ and \eqref{eq:twisting}, we have the following simple Deligne pairing interpretation of associated hypersurfaces:
\begin{equation}
\cA_{n+1-r}=\langle c_1(L)^r \cdot c_{n+1-r}(T_{X/S}^\vee\otimes L)\rangle_{X/S}.
\end{equation}
For $r=n$, this recovers the isomorphism part of the main result of \cite{kapadia}. Note that we also obtain an explicit expression for the degrees
\begin{equation}
d_{n+1-r}=\sum_{p=0}^{n+1-r}\binom{n+1-p}{r}\int_{X_s} c_p(T_{X_s}^\vee)\cdot c_1(L_s)^{n-p}.
\end{equation}
\end{example}

\begin{example}One also obtains a mixed version of the higher associated hypersurfaces, generalising the mixed Hurwitz forms of \cite{mixedhurwitz}, which give for example
\begin{equation}
\cA = \langle c_1(K_{X/S} + L_1 + \dots + L_n)\cdot c_1(L_1)\cdot\dots\cdot c_1(L_n)\rangle_{X/S}.
\end{equation}
Likewise, one can mix the previous two examples, involving multiple vector bundles in higher associated hypersurfaces. The statement of Theorem \ref{thm:main} more generally applies to any situation involving projective duality, and is thus 
\end{example}

\bibliographystyle{alpha}
\bibliography{bib}

\end{document}